\documentclass[12pt, reqno]{amsart}
\usepackage{fancyhdr}
\usepackage{a4wide, amsthm, amscd, amsfonts, amssymb, graphicx,tikz, color, environ}
\usepackage{amsaddr}
\usepackage{amsmath,amssymb,amsfonts,amsthm, graphicx,color,fancyhdr}
\usepackage[latin1]{inputenc}
\usepackage{enumerate}

\usepackage{psfrag}
\usepackage[colorlinks=true, pdfstartview=FitV, allcolors=black]{hyperref}
\usepackage{psfrag,caption}
\usepackage[all]{xy}
\usepackage{comment}
\usepackage{appendix}
\usepackage{verbatim}
\usepackage{mathrsfs}
\usepackage{pdfsync}
\usepackage{graphicx,epsfig}
\newtheorem{thm}{Theorem}[section]

\newtheorem{lem}[thm]{Lemma}
\theoremstyle{stylename}

\let\oldproofname=\proofname
\renewcommand{\proofname}{\rm\bf{\oldproofname}}

 \newtheorem{lemma}[thm]{Lemma}

\theoremstyle{definition}
 
 \newtheorem{rmk}[thm]{Remark}
 \newtheorem{defn}{Definition}[section]

\begin{document}

\title[ Upper bound for Neumann eigenvalue]
 {An upper bound for the first nonzero Neumann eigenvalue}
\author{Sheela Verma}
\address{Tata Institute of Fundamental Research \\ Centre For Applicable Mathematics \\ Bangalore, India}
\email{sheela.verma23@gmail.com}



\begin{abstract}
Let $\mathbb{M}$ denote a complete, simply connected Riemannian manifold with sectional curvature $K_{\mathbb{M}} \leq k$ and Ricci curvature $\text{Ric}_{\mathbb{M}} \geq (n-1)K$, where $k,K \in \mathbb{R}$. Then for a bounded domain $\Omega \subset\mathbb{M}$ with smooth boundary, we prove that  the first nonzero Neumann eigenvalue $\mu_{1}(\Omega) \leq \mathcal{C} \mu_{1}(B_{k}(R))$. Here $B_{k}(R)$ is a geodesic ball of radius $R > 0$ in the simply connected space form $\mathbb{M}_{k}$ such that vol$(\Omega)$ = vol$(B_{k}(R))$, and $\mathcal{C}$ is a constant which depends on the volume, diameter of $\Omega$ and the dimension of $\mathbb{M}$. 
\end{abstract}
\maketitle
\section{Introduction}

Let $M$ be a complete connected Riemannian manifold and $\Omega \subset {M}$ be a bounded domain with smooth boundary $\partial \Omega$. The Neumann eigenvalue problem on $\Omega$ is to find all real numbers $\mu(\Omega)$ for which there exists a nontrivial function $\varphi \in C^{2}(\Omega) \cap C^{1}(\bar{\Omega})$ such that
\begin{align}\label{Neumann problem}
\begin{array}{rcll}
\Delta \varphi &=&\mu \varphi & \text{ in } \Omega\\
\frac{\partial \varphi}{\partial \nu} &=& 0 & \text{ on } \partial \Omega
\end{array}
\end{align}
where $\nu$ is the outward unit normal to $\partial \Omega$. This problem has discrete and real spectrum $0 = \mu_{0}(\Omega) < \mu_{1}(\Omega) \leq \mu_{2}(\Omega) \leq \cdots \nearrow \infty$. In this article, we are interested in finding an upper bound for the first nonzero eigenvalue of problem \eqref{Neumann problem}. The variational characterization of $\mu_{1}(\Omega)$ is given by
\begin{align*}
\mu_{1}(\Omega) = \inf_{\varphi} \left\lbrace \frac{\int_{\Omega} \vert \nabla \varphi \vert ^{2} dV}{\int_{\Omega}  \varphi ^{2} dV} \mid \int_{\Omega}  \varphi \ dV = 0 \right\rbrace. 
\end{align*}

Among all simply connected planar domains of given area, Szeg\"{o} \cite{S} proved that the ball maximizes the first nonzero Neumann eigenvalue. This result was extended to arbitrary bounded domains in $\mathbb{R}^{n}$ by Weinberger \cite{W}. Using the idea of \cite{W}, Ashbaugh and Benguria \cite{AB} showed that Szeg\"{o}-Weinberger result \cite{S, W} also holds for bounded domains contained in a hemisphere of sphere $\mathbb{S}^{n}$. Later, Xu \cite{Xu} and Aithal-Santhanam \cite{AS} proved that the same result is also true for bounded domains in hyperbolic space $\mathbb{H}^{n}$ and rank-$1$ symmetric spaces, respectively.

Recently, Wang \cite{KW} considered problem \eqref{Neumann problem} on bounded domains in an $n$-dimensional complete, simply connected Riemannian manifold $\mathbb{M}$ whose sectional curvature is bounded from above by $k$ and Ricci curvature is bounded from below by $(n-1)K$ for $k,K \in \mathbb{R}$. Under certain assumptions on the size of domain $\Omega$, Wang \cite{KW} proved, 
\begin{align} \label{ineq:KWang}
\mu_{1}(\Omega) \leq \left( \frac{\sin_{K}(d)}{\sin_{k}(d)}\right)^{2n-2} \mu_{1}(\Omega^{*}),
\end{align}
where $\Omega^{*}$ denotes a geodesic  ball in a Riemannian manifold of constant curature $k$ such that vol$(\Omega)$ = vol$(\Omega^{*})$. Here $d$ is the diameter of $\Omega$ and the function $\sin_{m}(r)$ is defined as
\begin{align*}
\sin_{m}(r) :=
\begin{cases}
\frac{\sin \sqrt{m} r}{\sqrt{m}}, & m>0, \\
r, & m=0, \\
\frac{\sinh \sqrt{-m} r}{\sqrt{-m}}, & m<0. 
\end{cases}
\end{align*}
Basically, author in \cite{KW} used spherical symmetrization technique to prove the above result. \\

In this article, we give an upper bound for the first nonzero Neumann eigenvalue. The main contribution and properties of this work are summarised as follows:
\begin{enumerate} [(i)]
\item Under the same assumptions on $\Omega$ as in \cite{KW}, we prove that $\mu_{1}(\Omega) \leq \mathcal{C} \ \mu_{1}(B_{k}(R))$, where the constant $\mathcal{C}$(defined in Section $3$) depends on $k, K, n, d$, volume of $\Omega$, and $B_{k}(R)$ is a geodesic ball of radius $R$ in the space form of curvature $k$ such that vol$(\Omega)$ = vol$(B_{k}(R))$.
\item For a domain of given volume and arbitrary large diameter, the constant $\mathcal{C}$ in our result is smaller than the constant $\left( \frac{\sin_{K}(d)}{\sin_{k}(d)}\right)^{2n-2}$ appears in \cite{KW} (see Remark \ref{rmk:comparison}).
\item We also observe that arguments given in \cite{KW} are analytical in nature however our arguments are more geometrical.
\end{enumerate}

\section{Preliminaries}
In this section, we recall the notion of center of mass and some basic results which are important to prove our result. 

Let $(M, g)$ be an $n$-dimensional complete Riemannian manifold. For a point $p \in M$, we denote the
convexity radius of $M$ at $p$ by $c(p)$. For a domain $\Omega \subset B(p, c(p))$, we denote the convex hull of $\Omega$ by hull$(\Omega)$. Let $\exp_{p} : T_{p}M \rightarrow M$ be the exponential map. The following lemma gives the existence of a center of mass of a bounded domain in $M$.

\begin{lem} \label{lem: centre of mass}
Let $\Omega$ be a bounded domain in $(M,g)$ and is contained in $B(q,c(q))$ for some $q \in M$. Let $G$ be a continuous function on $\left[0, 2 c(q) \right]$ such that $G$ is also positive on $\left(0, 2 c(q) \right) $. Then there exists a point $p \in \text{hull}(\Omega)$ such that 
\begin{align*}
\int_{\Omega} G(r_{p}(x)) \frac{\exp_{p}^{-1}(x)}{r_{p}(x)} dV = 0,
\end{align*}
where $r_{p}(x)$ denotes the distance between $p$ and $x$ in $M$.
\end{lem}

For a proof see \cite{AS}.
\begin{defn}
The point $p$ in the above lemma is called a center of mass of domain $\Omega$ with respect to the function $G$.
\end{defn}

Next we fix some notations which we use throughout the manuscript.

Let $\mathbb{M}$ be a complete, simply connected Riemannian manifold of dimension $n$ with sectional curvature $K_{\mathbb{M}} \leq k$ and Ricci curvature bounded from below by $(n-1)K$. We denote the simply connected space form of constant curvature $s$ by $\mathbb{M}_{s}$. Let $\Omega$ be a bounded domain in $\mathbb{M}$ with smooth boundary and vol$(\Omega)$ = vol$(B_{k}(R))$ = vol$(B_{K}(R'))$, where $B_{k}(R)$ and $B_{K}(R')$ denote balls of radius $R$ (in $\mathbb{M}_{k}$) and $R'$ (in $\mathbb{M}_{K}$), respectively. We further assume that balls $B_{K}(R')$ and $B_{k}(R)$ are centered at $p_{K} \in \mathbb{M}_{K}$ and $p_{k} \in \mathbb{M}_{k}$, respectively. For $y \in \mathbb{M}_{k}$, $r_{p_{k}}(y)$ represents the distance between $p_{k}$ and $y$ in $\mathbb{M}_{k}$. Similarly, $r_{p_{K}}(y')$ is the distance between $p_{K}$ and $y'$ in $\mathbb{M}_{K}$, where $y' \in \mathbb{M}_{K}$. Observe that $R' \leq R$. For $k>0$, we impose the following condition on $\Omega$.
\begin{enumerate}[(A)]
\item $\begin{aligned} [t]
d = \text{diam}(\Omega) = \text{diam}(\text{hull} (\Omega)) < \min \left\lbrace \frac{\pi}{2 \sqrt{k}}, \text{ injectivity radius of } \mathbb{M}\right\rbrace. 
\end{aligned}$

\item $\begin{aligned} [t]
\text{vol}_{\mathbb{M}}(\text{hull}(\Omega)) \leq \frac{\text{vol}_{\mathbb{M}_{k}}(\mathbb{M}_{k})}{2} . 
\end{aligned}$
\end{enumerate}

\begin{rmk}
For $k>0$, assumption $(A)$ assures that exponential map $\exp_{p}$ is a diffeomorphism onto hull$(\Omega)$ for all $p \in \text{hull}(\Omega)$, and we require condition $(B)$ so that we can ultimately work in the hemisphere of $\mathbb{M}_{k}$.
\end{rmk}
Now we state some properties of $\mu_{1}(B_{k}(R))$, the first nonzero Neumann eigenvalue of $B_{k}(R)$. 

Recall that by separation of variable technique, the first nonzero Neumann eigenvalue of $B_{k}(R)$ is the first eigenvalue of
\begin{align} \label{eqn:separationofvariable}
\begin{array} {rcl}
- F''(r) - \frac{(n-1) \sin_{k}'(r)}{\sin_{k}(r)} F'(r) + \frac{n-1}{\sin_{k}^{2}(r)} F(r) = \mu F(r) \\
F(0)=0, \quad F'(R) = 0.
\end{array}
\end{align}

Let $f(r)$ be the eigenfunction of \eqref{eqn:separationofvariable} corresponding to $\mu = \mu_{1}(B_{k}(R))$. Then $f(r(y)) \frac{y_{i}}{r(y)}$, $1 \leq i \leq n$ is an eigenfunction corresponding to $\mu_{1}(B_{k}(R))$, where $(y_{1}, y_{2}, \ldots, y_{n})$ is the geodesic polar coordinates of $y \in B_{k}(R)$ with respect to the center of $B_{k}(R)$ and $r(y)$ represents the distance between $y$ and the center of $B_{k}(R)$. The function $f(r)$ satisfies the following properties.
\begin{lemma}
$f(r)$ is an increasing function of $r$ and $\left( f'(r)\right)^{2} + \frac{n-1}{\sin_{k}^{2}(r)} f^{2}(r)$ is a decreasing function of $r$.
\end{lemma}
For more details see \cite{AB, Xu}.

The following lemmas are useful in proving our main result.
\begin{lemma}
 For a fix point $q \in \mathbb{M}$, let $(x_1, x_2, \ldots, x_n)$ denote the geodesic polar coordinates with respect to $q$. Let $S(r)$ be the geodesic sphere of radius $r$ centered at $q$. For $k>0$, we further assume
\begin{align*}
r < \min \left\lbrace \frac{\pi}{2 \sqrt{k}}, \text{ injectivity radius of } \mathbb{M}\right\rbrace. 
\end{align*}
Then
\begin{align}
\sum_{i=1}^{n} \left| \nabla^{S(r)} \left( \frac{x_{i}}{r} \right)\right|^{2} \leq \frac{n-1}{\sin_{k}^{2}(r)}.
\end{align}
\end{lemma}
For more general statement and details see \cite{BS}.

\begin{lemma} \label{sin inequality}
For $K, k$ and $\sin_{m}(r)$ defined as above, the function $\frac{\sin_{K}(r)}{\sin_{k}(r)}$ is an increasing function of $r$.
\end{lemma}

Using the exponential map, next we construct a domain $\Omega_{K}$ in $\mathbb{M}_{K}$ from $\Omega$. This construction is further used in the proof of our main result.

For $p \in \mathbb{M}$ , let $W \subset T_{p}(\mathbb{M})$ such that $\Omega = \exp_{p}(W)$. Fix an isometry $i : T_{p}(\mathbb{M}) \longrightarrow T_{p_{K}}(\mathbb{M}_{K})$ for $p_{K} \in \mathbb{M}_{K}$ and denote $\Omega_{K} = \exp_{p_{K}}(i(W))$. Note that corresponding to each $q \in \Omega$, there exists $\bar{q} \in \Omega_{K}$ such that $q = \exp_{p} u$ and $\bar{q} = \exp_{p_{K}}(i(u))$ and vice versa. Thus in terms of geodesic polar coordinates, $\Omega$ can be written as
\begin{align*}
\Omega &= \left\lbrace \left( r, u\right): u \in T_{p}\mathbb{M},\Vert u \Vert = 1 , r \in (r_{1}(u), r_{2}(u))\cup (r_{3}(u), r_{4}(u)) \cup \cdots \right. \\& \quad\left. \cup (r_{m(u)-1}(u), r_{m(u)}(u)) \right\rbrace 
\end{align*}
Similarly
\begin{align*}
\Omega_{K} & = \left\lbrace \left( r, \bar{u}\right): \bar{u} \in T_{p_{K}}\mathbb{M}_{K}, \Vert u_{K} \Vert = 1 , r \in (r_{1}(\bar{u}), r_{2}(\bar{u})) \cup (r_{3}(\bar{u}), r_{4}(\bar{u})) \cup \cdots \right. \\ & \quad \left.  \cup (r_{m(\bar{u})-1}(\bar{u}), r_{m(\bar{u})}(\bar{u})) \right\rbrace \\
& = \left\lbrace \left( r, i(u)\right): u \in T_{p}\mathbb{M}, \Vert u \Vert = 1 , r \in (r_{1}(u), r_{2}(u)) \cup (r_{3}(u), r_{4}(u)) \cup \cdots \right. \\  & \quad \left. \cup (r_{m(u)-1}(u), r_{m(u)}(u)) \right\rbrace.
\end{align*}

Denote $I_{u} = I_{\bar{u}}= (r_{1}, r_{2}) \cup (r_{3}, r_{4}) \cup \cdots \cup (r_{m(u)-1}, r_{m(u)})$. Let $\phi$ and $\phi_{K} (= \sin_{K}^{n-1}(r))$ be the volume density functions of $\mathbb{M}$ and $\mathbb{M}_{K}$ along the radial geodesics starting from $p$ and $p_{K}$,  respectively. Observe that by Gunther volume comparison theorem,  vol$(\Omega) = \text{vol} (B_{K}(R')) \leq \text{vol} (\Omega_{K})$.
\section{Statement and proof of the main result}
\begin{thm}
With all assumptions on $\Omega$ given in Section $2$, the first nonzero Neumann eigenvalue on $\Omega$ satisfies the following inequality
\begin{align*}
\mu_{1}(\Omega) \leq \mathcal{C} \  \mu_{1}(B_{k}(R)),
\end{align*}
where
\begin{align*}
\mathcal{C} = \left( \frac{\sin_{K}(R)}{\sin_{k}(R)}\right)^{n-1} \left( \frac{\sin_{K}(d)}{\sin_{k}(d)}\right)^{n-1}  \frac{\int_{ B_{k}(R)} f^{2}(r_{p_{k}}) dV}{\int_{ B_{K}(R')} f^{2}(r_{p_{K}}) dV}.
\end{align*}
Further, if $k = K$ then the constant $\mathcal{C}$ is equal to $1$, and the above bound is sharp.
\end{thm}

\begin{rmk} \label{rmk:comparison}
Since the function $\frac{\sin_{K}(r)}{\sin_{k}(r)}$ is increasing and unbounded, for $k <0$ and among all domains of given volume, the constant factor $\left( \frac{\sin_{K}(d)}{\sin_{k}(d)}\right)^{2n-2}$ in \eqref{ineq:KWang} is larger than the constant $\mathcal{C}$ for domains of arbitrary large diameter.
\end{rmk}

\begin{proof}
Define
\begin{align*}
h(r) :=
\begin{cases}
f(r), \quad r \leq R \\
f(R), \quad r \geq R.
\end{cases}
\end{align*}
Note that the function $h(r)$ is continuous and positive function on $[0, \infty) $. Let $p \in \text{hull} (\Omega)$ be a centre of mass with respect to the function $h(r)$ and $(x_{1}, x_{2},\ldots, x_{n})$ denote the geodesic polar coordinates of $x \in \Omega$ with respect to the point $p$. For our convenience, in this proof we denote $r_{p}(x)$ (given in Lemma \ref{lem: centre of mass}) by $r$, $r_{p_{k}}(y)$ by $r_{p_{k}}$  and $r_{p_{K}}(y')$ by $r_{p_{K}}$. Then for all $1 \leq i \leq n$,
\begin{align*}
\int_{\Omega} h(r) \frac{x_{i}}{r} dV = 0.
\end{align*}
Now by the variational characterization of $\mu_{1}(\Omega)$, we have
\begin{align*}
\mu_{1}(\Omega) \sum_{i=1}^{n} \int_{\Omega} \left( h(r) \frac{x_{i}}{r}\right)^{2} dV \leq \sum_{i=1}^{n} \int_{\Omega} \left| \nabla \left( h(r) \frac{x_{i}}{r} \right)\right|^{2} dV, \\
\mu_{1}(\Omega) \int_{\Omega} h^{2}(r) dV \leq \sum_{i=1}^{n} \int_{\Omega} \left| \nabla \left( h(r) \frac{x_{i}}{r} \right)\right|^{2} dV.
\end{align*}

Next we find an estimate for $\int_{\Omega} h^{2}(r) dV$.
\begin{align*}
\int_{\Omega} h^{2}(r) dV &= \int_{U_{p}\mathbb{M}} \int_{I_{u}} h^{2}(r) \ \phi(r, u) \ dr \ du \\
& \geq \int_{U_{p}\mathbb{M}} \int_{I_{u}} h^{2}(r) \ \sin_{k}^{n-1}(r) \ dr \ du \\
& = \int_{U_{p}\mathbb{M}} \int_{I_{u}} h^{2}(r) \   \frac{\sin_{k}^{n-1}(r)}{\sin_{K}^{n-1}(r)} \sin_{K}^{n-1}(r) \ dr \ du.
\end{align*}
Since $\frac{\sin_{k}(r)}{\sin_{K}(r)}$ is a decreasing function of $r$,
\begin{align*} 
\int_{\Omega} h^{2}(r) dV & \geq \frac{\sin_{k}^{n-1}(d)}{\sin_{K}^{n-1}(d)} \int_{U_{p}\mathbb{M}} \int_{I_{u}} h^{2}(r) \  \sin_{K}^{n-1}(r) \ dr \ du \\ 
& = \frac{\sin_{k}^{n-1}(d)}{\sin_{K}^{n-1}(d)} \int_{U_{p_{K}}\mathbb{M}_{K}} \int_{I_{\bar{u}}} h^{2}(r) \  \sin_{K}^{n-1}(r) \ dr \ d\bar{u} . 
\end{align*}

\begin{align*}
 \int_{\Omega} h^{2}(r) dV & \geq \frac{\sin_{k}^{n-1}(d)}{\sin_{K}^{n-1}(d)} \int_{\Omega_{K}} h^{2}(r_{p_{K}}) \ dV \\ 
& = \frac{\sin_{k}^{n-1}(d)}{\sin_{K}^{n-1}(d)} \left( \int_{\Omega_{K} \cap B_{K}(R')} h^{2}(r_{p_{K}}) dV +\int_{\Omega_{K} \setminus \left( \Omega_{K} \cap B_{K}(R')\right) } h^{2}(r_{p_{K}}) \ dV \right). 
\end{align*}

Using the fact that $h(r)$ is an increasing function of $r$, we obtain
\begin{align} \nonumber
\int_{\Omega} h^{2}(r) & \geq \frac{\sin_{k}^{n-1}(d)}{\sin_{K}^{n-1}(d)} \left( \int_{ B_{K}(R')} h^{2}(r_{p_{K}}) dV - \int_{B_{K}(R') \setminus \left( \Omega_{K} \cap B_{K}(R')\right) } h^{2}(R') \ dV \right.   \\ \nonumber
& \qquad \left.  +\int_{\Omega_{K} \setminus \left( \Omega_{K} \cap B_{K}(R')\right) } h^{2}(R') \ dV \right) \\ \label{ineq:fun}
& \geq \frac{\sin_{k}^{n-1}(d)}{\sin_{K}^{n-1}(d)} \int_{ B_{K}(R')} h^{2}(r_{p_{K}}) dV.
\end{align}
The last inequality follows from the fact that $\text{vol} (B_{K}(R')) \leq \text{vol} (\Omega_{K})$.\\

Now we obtain an upper bound for $\sum_{i=1}^{n} \int_{\Omega} \left| \nabla \left( h(r) \frac{x_{i}}{r} \right)\right|^{2} dV $.

\begin{align*}
\sum_{i=1}^{n} \int_{\Omega} \left| \nabla \left( h(r) \frac{x_{i}}{r} \right)\right|^{2} dV & = \int_{\Omega} \left( h^{2}(r)\sum_{i=1}^{n} \left| \nabla^{S(r)} \left( \frac{x_{i}}{r} \right)\right|^{2} + \left( h'(r)\right) ^{2} \right) dV \\
& \leq \int_{\Omega} \left( \frac{n-1}{\sin_{k}^{2}(r)} h^{2}(r)+ \left( h'(r)\right) ^{2} \right) dV.
\end{align*}
Denote $G(r)= \frac{n-1}{\sin_{k}^{2}(r)} h^{2}(r)+ \left( h'(r)\right) ^{2}$. Since $G(r)$ is a decreasing function of $r$,
\begin{align*} 
\sum_{i=1}^{n} \int_{\Omega} \left| \nabla \left( h(r) \frac{x_{i}}{r} \right)\right|^{2} dV & \leq \int_{\Omega \cap B(R)} G(r) dV +\int_{\Omega \setminus \left( \Omega \cap B(R)\right) } G(r) \ dV \\ 
& = \int_{B(R)} G(r) dV - \int_{B(R) \setminus \left( \Omega \cap B(R)\right) } G(r) \ dV +\int_{\Omega \setminus \left( \Omega \cap B(R)\right) } G(r) \ dV \\
 & \leq \int_{B(R)} G(r) dV \\ 
& = \int_{U_{p}\mathbb{M}} \int_{0}^{R} G(r) \ \phi(r, u) \ dr \ du.
\end{align*}
Using Gunther volume comparison theorem to get
\begin{align} \nonumber
\sum_{i=1}^{n} \int_{\Omega} \left| \nabla \left( h(r) \frac{x_{i}}{r} \right)\right|^{2} dV & \leq \int_{U_{p}\mathbb{M}} \int_{0}^{R} G(r) \ \sin_{K}^{n-1}(r) \ dr \ du \\ \nonumber
& \leq \left( \frac{\sin_{K}(R)}{\sin_{k}(R)}\right)^{n-1}  \int_{U_{p_{k}}\mathbb{M}_{k}} \int_{0}^{R} G(r) \ \sin_{k}^{n-1}(r) \ dr \ d\bar{u} \\ \label{ineq:grad}
& = \left( \frac{\sin_{K}(R)}{\sin_{k}(R)}\right)^{n-1}  \int_{B_{k}(R)} G(r_{p_{k}}) dV.
\end{align}
We have used Lemma \ref{sin inequality} to conclude the second last inequality. By combining \eqref{ineq:fun} and \eqref{ineq:grad}, we obtain
\begin{align*}
\mu_{1}(\Omega) & \leq \left( \frac{\sin_{K}(R)}{\sin_{k}(R)}\right)^{n-1} \left( \frac{\sin_{K}(d)}{\sin_{k}(d)}\right)^{n-1} \frac{\int_{B_{k}(R)} \left( \frac{n-1}{s_{k}^{2}(r_{p_{k}})} h^{2}(r_{p_{k}})+ \left( h'(r_{p_{k}})\right) ^{2} \right) dV}{\int_{ B_{K}(R')} h^{2}(r_{p_{K}}) dV}\\
& = \mathcal{C} \ \frac{\int_{B_{k}(R)} \left( \frac{n-1}{s_{k}^{2}(r_{p_{k}})} f^{2}(r_{p_{k}})+ \left( f'(r_{p_{k}})\right) ^{2} \right) dV}{\int_{ B_{k}(R)} f^{2}(r_{p_{k}}) dV}\\
& = \mathcal{C}  \ \mu_{1}(B_{k}(R)),
\end{align*}
where
\begin{align*}
 \mathcal{C}  = \left( \frac{\sin_{K}(R)}{\sin_{k}(R)}\right)^{n-1} \left( \frac{\sin_{K}(d)}{\sin_{k}(d)}\right)^{n-1}  \frac{\int_{ B_{k}(R)} f^{2}(r_{p_{k}}) dV}{\int_{ B_{K}(R')} f^{2}(r_{p_{K}}) dV}.
\end{align*}
This completes the proof of the theorem.
\end{proof}


\begin{thebibliography}{9}
\bibitem{AS} A. R. Aithal, G. Santhanam, Sharp upper bound for the first non-zero Neumann eigenvalue for bounded domains in rank-1 symmetric spaces, \emph{Transactions of the American Mathematical Society} \textbf{348}(10) 3955-3965 (1996).

\bibitem{AB} M.S. Ashbaugh, R.D. Benguria, Sharp upper bound to the first nonzero Neumann eigenvalue for bounded domains in spaces of constant curvature, \emph{Journal of the London Mathematical Society} \textbf{52}(2) 402-416 (1995).

\bibitem{BS} Binoy, G. Santhanam, Sharp upperbound and a comparison theorem for the first nonzero Steklov eigenvalue, \emph{Journal of Ramanujan Mathematical Society} \textbf{29}(2) 133-154 (2014).

\bibitem{S} G. Szeg\"{o}, Inequalities for certain eigenvalues of a membrane of given area, \emph{Journal of Rational Mechanics and Analysis} \textbf{3} 343-356 (1954).

\bibitem{KW} K. Wang, An upper bound for the second Neumann eigenvalue on Riemannian manifolds, \emph{Geometriae Dedicata} \textbf{201}(1) 317-323 (2019).

\bibitem{W} H. F. Weinberger, An isoperimetric inequality for the N-dimensional free membrane problem, \emph{Journal of Rational Mechanics and Analysis} \textbf{5}(4) 633-636 (1956).

\bibitem{Xu} Y. Xu, The first nonzero eigenvalue of Neumann problem on Riemannian manifolds, \emph{The Journal of Geometric Analysis } \textbf{5}(1) 151-165 (1995).


\end{thebibliography}
\end{document}